\documentclass[12pt]{amsart}
\usepackage{amssymb}
\usepackage{amscd}
\usepackage{amsmath}
\usepackage{amsthm}
\usepackage[nobysame,initials]{amsrefs}

\newcommand{\IGNORE}[1] {}

\newcommand{\Q}{\mathbb Q}
\newcommand{\Z}{\mathbb Z}

\newcommand{\C}{\mathbb C}

\newcommand{\A}{\mathbb A}

\newcommand{\N}{\mathbb N}

\newcommand{\foralmostall}{{\text for \, all \, but \, finitely \, many}}

\makeatletter
\def\l@section{\@tocline{1}{4pt}{1pc}{}{}}
\def\l@subsection{\@tocline{2}{0pt}{2pc}{5pc}{}}
\makeatother
\makeatletter
\def\imod#1{\allowbreak\mkern10mu({\operator@font mod}\,\,#1)}
\makeatother

\title{A mild Tchebotarev theorem for GL$(n)$}
\author{Dinakar Ramakrishnan}\thanks{Partly supported by the NSF grant DMS-0701089}
\begin{document}
\maketitle

\medskip

\begin{flushright}
{\it In memory of Steve Rallis}
\end{flushright}


\section*{Introduction}

\medskip

As it is well known, the Tchebotarev density theorem implies that
two irreducible $\ell$-adic representations $\rho_\ell$,
$\rho'_\ell$ of the absolute Galois group of a number field $K$ are
isomorphic if the corresponding characteristic polynomials of
Frobenius elements agree on a set $S$ of primes of
density $1$. It is then natural to ask, in view of the Langlands
conjectures, whether an analogous assertion holds for cuspidal
automorphic representations of GL$_n(\A_K)$. The object of this Note
is to establish such an automorphic analogue for a simple, but useful, class
of $S$ of density $1$. To be precise, we prove the following:

\medskip

\noindent{\bf Theorem A} \, \it Let $K/k$ be a cyclic extension
of number fields of degree
a prime $p$, and let $\Sigma^1_{K/k}$ denote the set
of primes $v$ of $K$ which are of degree $1$ over $k$. Suppose
$\pi$, $\pi'$ are cusp forms on GL$(n)/K$ such that $\pi_v \, \simeq
\pi'_v$, for all but a finite number of $v$ in $\Sigma^1_{K/k}$.
Then $\pi, \pi'$ are twist equivalent. More precisely,
they have isomorphic base changes over the cyclotomic extension
$K(\zeta)$, where $\zeta$ is a non-trivial $p$-th root of unity. \rm

\medskip

We refer to the book \cite{AC} for facts on solvable base change for GL$(n)$ due to Arthur and Clozel.

When we say that $\pi, \pi'$ are twist equivalent, we mean $\pi' \simeq \pi\otimes \chi$ for a
finite order character $\chi$ of (the idele classes of) $K$. In particular, if $n$ is
relatively prime to $p-1$, or if the
conductors of $\pi, \pi'$ are prime to $p$, we may conclude
even that
$\pi, \pi'$ are isomorphic (over $K$). When $p=2$, we thus get the following:

\medskip

\noindent{\bf Corollary B} \, \it Let $K/k$ be a quadratic extension of number fields.
Then any cuspidal automorphic representation $\pi$ of GL$_n(\A_K)$ is determined (up to isomorphism)
by its components $\pi_v$ for all (but a finite number of) places $v$ of degree $1$ over $k$.\rm

\medskip

Clearly, Theorem A refines the strong
multiplicity one theorem, which gives the desired global isomorphism if $\pi_v \simeq \pi'_v$
for all but a {\it finite} number of $v$.
(\cite{JS}). For GL$(2)$, there is a stronger result known,
requiring the isomorphism $\pi_v \simeq \pi'_v$ only for a set $S'$
of $v$ of density $> 7/8$ (\cite{Ra1}). For GL$(n)$ with $n>2$, we conjectured elsewhere that such a stronger
result should hold with $7/8$ replaced by $1- 1/2n^2$, which is a theorem for $\pi$ attached to an
$\ell$-adic representation $\rho_\ell$ by an elegant result of Rajan (\cite{Raj1}). We are {\it far} from
such a precise result for those cusp forms $\pi$ on GL$(n)$, $n \geq 3$, which are not known to be associated
to such a $\rho_\ell$.

\medskip

Given a finite cyclic extension $K/k$, if $G$, resp. $\tilde G$, is
a reductive group over $k$, resp. $K$, such that $\tilde G =
G\times_k K$, let us say that a cuspidal automorphic representation
$\pi$ of $G(\A_k)$ admits a {\it soft base change} to $K$ if there is
an automorphic representation $\Pi$ of $\tilde G(\A_K)$ such that
for all but a finite number of primes $v$ in $\Sigma_{K/k}^1$, we have
$\Pi_v \simeq \pi_u$, where
$u$ is the prime of $k$ below $v$. When $\tilde G$ is GL$(n)/K$,
Theorem A says that a soft base change $\Pi$ is unique up to twisting by a character of
$K$ which is trivial when pulled back by norm to the $p$-cyclotomic extension of $K$; in particular, $\Pi$ is determined up to isomorphism for $p=2$
when cuspidal. An initial impetus for it came from a question
asked independently by J.~Getz and D.~Whitehouse. Now Theorem A has been used
(for quadratic extensions) by B.~Feigon, K.~Martin and D.~Whitehouse in their paper (\cite{FMW})
on the periods and non-vanishing of $L$-functions of GL$(2n)$, and 
in Wei Zhang's work on the Gross-Prasad conjecture (\cite{Zh}).

\medskip

Now a few words about the proof of Theorem A. A well known, basic
theorem of Luo, Rudnick and Sarnak (\cite{LRS}), which is of importance to us, says that for any
cusp form $\pi$ on GL$(n)/K$, the coefficient $a_v$ of $\pi$ at any
unramified $v$ satisfies the bound $\vert a_v\vert <
(Nv)^{1/2-1/(n^2+1)}$. (What is essential for us is that $a_v$ is bounded in absolute value by
$(Nv)^{1/2-t_n}$ for a fixed positive number $t_n$ independent of $v$, not the exact shape of $t_n$.) As it has been noted and used already by Rajan (\cite{Raj2}), feeding this into the known analytic framework, it suffices, under our hypotheses, to prove that for all
but a finite number of $v$ whose degree lies in $[2,(n^2+1)/2]$,
$\pi_v$ and $\pi'_v$ are isomorphic. We cannot achieve this
directly, but can show, using some Kummer theory, that it holds for the base changes $\pi_L, \pi'_L$
to a carefully chosen
solvable extension $L$ of $K'=K(\zeta)$, which will be a compositum (over $K$) of a finite number of disjoint $p^r$-extensions $L^{(1)}, L^{(2)}, \dots$ with $2p^r>n^2+1$; each $L^{(j)}$ will be a nested chain of cyclic $p^2$-extensions (see section 4). From this data we prove by descent that
$\pi_{K'}$ and $\pi'_{K'}$ are isomorphic. There is an added subtlety if $\pi_{K'}$ or $\pi'_{K'}$ is not cuspidal, and this forces us to work with isobaric sums of
unitary cuspidal automorphic representations, which are analogues of semisimple Galois representations of pure weight. These steps together form the core of the argument. It should be stressed that since the basic analytic method is by now standard, given Rajan's work (\cite{Raj2}) making use of \cite{LRS}, what is new here is the use of base change to a suitable chain of $p$-power extensions to achieve the requisite isomorphism, followed by careful descent.

In another paper (\cite{Ra2}), we extend Theorem A non-trivially to the case of an arbitrary Galois extension $K/k$.  The main idea there is quite different and replaces explicit Kummer theory with a fuller use of class field theory, in particular the Tate cohomology and duality. We hope that it is still of interest to have just the cyclic case published, at least because the proof is simpler and more accessible.

\medskip

We thank the many people who have shown interest in this work over the past few years, especially to those who have used it and have encouraged, like K.~Martin, to have it published. Thanks are also due
to the NSF for partial support through the grants DMS-0701089 and DMS-1001916. This article is dedicated to the memory of Steve Rallis from whom this author learnt a lot in conversations over the years.

\bigskip

\section{Basic Facts: A Review}

\medskip

Let $F$ be a global field with ad\`ele ring $\A_F$.  Let $\Sigma_F$
denote the set of all
places of $F$.  If $v\in\Sigma_F$ is finite, let $q_v$ denote the cardinality
of the residue field at $v$.  For $n \geq 1$, let $A_0(n,F)$ denote the set of
isomorphism classes of irreducible unitary, cuspidal automorphic
representations of GL$(n,\A_F)$.  Every $\pi$ representing a class
in $A_0(n,F)$ is (isomorphic to) a tensor product $\otimes_v,
\pi_v$, where $v$ runs over all the places of $F$, such that each
$\pi_v$ is an irreducible generic representation of GL$(n,F_v )$
and such that $\pi_v$ is unramified at almost all $v$.
The strong multiplicity one theorem (\cite{JS}) asserts that, for any {\it
finite} subset $S$ of $\Sigma_F$, $\pi$ is determined
up to isomorphism by the collection $\{\pi_v \, | \, v\not\in S\}$.

For any irreducible, automorphic representation $\pi$ of
$GL(n,\A_F),$ let $L(s, \pi) = L(s, \pi_{\infty})L(s, \pi_f)$ denote
the associated {\it standard} $L-$function of $\pi;$ it
has an Euler product expansion
$$
L(s,\pi) \, = \, \prod_v \, L(s, \pi_v),
$$
convergent in a right-half plane. If $v$ is a finite place where $\pi_v$ is
unramified, there is a corresponding semisimple (Langlands)
conjugacy class $A_v(\pi)$ (or $A(\pi_v)$) in GL$(n,\C)$ such that
$$
L(s,\pi_v) \, = \, {\rm {det}}(1-A_v(\pi)T)^{-1}|_{T=q_v^{-s}}.
$$
One may find a diagonal representative diag$(\alpha_{1,v}(\pi), ...
, \alpha_{n,v}(\pi))$ for $A_v(\pi),$ which is unique up to
permutation of the diagonal entries. Let $[\alpha_{1,v}(\pi), ...
, \alpha_{n,v}(\pi) ]$ denote the resulting unordered $n-$tuple.
One knows (by Godement-Jacquet) that for any non-trivial cuspidal
representation $\pi$ of GL$(n,\A_F)$,
$L(s,\pi)$ is entire.

\medskip

By Langlands's theory of Eisenstein series, one has a sum operation
$\boxplus$, called the isobaric sum (\cite{JS}): Given any
$m-$tuple of cuspidal representations $\pi_1, ..., \pi_m$ of
GL$(n_1,\A_F), ... ,$ GL$(n_m,\A_F)$ respectively, there exists an
irreducible, automorphic representation $\pi_1 \boxplus ...
\boxplus \pi_m$ of GL$(n,\A_F),$ $n \, = \, n_1 + ... + n_m$, which
is unique up to equivalence, such that for any finite set $S$ of
places,
$$
L^S(s, \boxplus_{j=1}^m \pi_j) \, = \, \prod_{j=1}^m L^S(s,
\pi_j).
$$
Call such a (Langlands) sum $\pi \simeq \boxplus_{j=1}^m \pi_j$,
with each $\pi_j$ cuspidal, an {\it isobaric} representation.

Denote by $\mathcal  A(n,F)$ the set, up to equivalence, of isobaric
automorphic representations of GL$_n(\A_F)$, and by
$\mathcal  A_u(n,F)$ the subset of isobaric sums of {\it unitary} cuspidal
automorphic representations. If $\pi=\boxplus_{i=1}^m \pi_i$, resp.
$\pi'=\boxplus_{j=1}^r \pi_j'$, is in
$\mathcal A_u(n,F)$, resp. $\mathcal A_u(n',F)$, with $\pi_i, \pi'_j$ unitary cuspidal,
we will need to consider the associated Rankin-Selberg $L$-function
$$
L(s, \pi \times \pi') \, = \, \prod_{i, j} \, L(s, \pi_i\times \pi'_j),
$$
with
$$
L(s, \pi_{i,v} \times \pi'_{j,v}) \, = \, {\rm {det}}(1-A_v(\pi_i)\otimes A_v(\pi'_j)T)^{-1}|_{T=q_v^{-s}}.
$$

\medskip

If $L(s)=\prod_{v\in\sum_\infty\cup\sum_f}L_v (s)$ is any global
$L$-function and $Y$ a set of places of $F$, then we will denote by
$L^Y(s)$ (resp. $L_Y(s)$) the product of $L_v (s)$ over all $v$
outside $Y$ (resp. in $Y$). We have the following basic result (\cite{JS}):

\medskip

\noindent{\bf Theorem 1.1} (Jacquet--Piatetski-Shapiro--Shalika,
Shahidi) \, \it Let $\pi=\boxplus_{i=1}^m \pi_i$,  $\pi'=\boxplus_{j=1}^r \pi_j'$
be in $\mathcal A_u(n,F)$, with $\pi_i, \pi'_j$ unitary cuspidal. Suppose $Y$ is a finite set
of places of $F$ containing the archimedean places
such that $\pi, \pi'$ are unramified
outside $Y$. Then $L^S(s, \pi \times \overline\pi')$ has a pole at $s=1$ iff for some $(i, j)$, $\pi_i$ is isomorphic to $\pi'_j$, in which case the pole of the factor $L(s, \pi_i \times \pi_j)$ is simple.\rm

\medskip

Here $\overline\pi'$ denotes the complex conjugate representation of $\pi'$, which, by unitarity, is the contragredient of $\pi'$.

\medskip

The general Ramanujan conjecture predicts that for any $\pi\in \mathcal A_u(n, F)$, $\pi_v$ is tempered at all $v$. In particular, if $v$ is a finite place where $\pi$ is unramified, the unordered $n$-tuple $\{\alpha_{1,v}(\pi), ...
, \alpha_{n,v}(\pi)\}$ representing $A_v(\pi)$ should satisfy $\vert \alpha_{i,v}\vert =1$ for every $i$. This is far from being proved, and the best known  bound to date (for general $n$) is given by the following:

\medskip

\noindent{\bf Theorem 1.2} (Luo--Rudnick--Sarnak \cite{LRS}) \, \it Let $\pi \in \mathcal A_u(n, F)$, and $v$ a finite place where $\pi$ is unramified, with $A_v(\pi)=\{\alpha_{1,v}(\pi), ...
, \alpha_{n,v}(\pi)\}$. Then for every $j\leq n$, one has
$$
\vert \alpha_{j,v}\vert \, < \, q_v^{\frac12-\frac{1}{n^2+1}}.
$$
\rm

\medskip

To be precise, Luo, Rudnick and Sarnak only address the case of cusp forms. But for $\pi \in \mathcal A_u(n,F)$, any $\alpha_j(\pi)$ must be associated to a cuspidal isobaric constituent $\pi_i$ on GL$(n_i)/F$ with $n_i \leq n$, and so the assertion above follows immediately from \cite{LRS}.

\medskip

We will also need the following (weak) version of the base change theorem for GL$(n)$:

\medskip

\noindent{\bf Theorem 1.3} (Arthur--Clozel \cite{AC}) \, \it Let $M/F$ be a finite extension of number fields obtained as a succession of cyclic extensions. Then for every $\pi \in \mathcal A_u(n,F)$, there exists a corresponding $\pi_M \in \mathcal A_u(n,M)$ such that for every finite place $v$ of $F$ where $\pi$ and $M$ are unramified, and for all places $w$ of $M$ dividing $v$, we have
$$
A_v(\pi)=\{\alpha_{1,v}, ...
, \alpha_{n,v}\} \, \implies \, A_w(\pi_{M})=\{\alpha_{1,v}^{f_v}, ...
, \alpha_{n,v}^{f_v}\},
$$
where $f_v=[M_w:F_v]$.
\rm

\medskip

A word of explanation may be helpful. In \cite{AC}, it is proved that for every cuspidal $\pi$, the base change $\pi_M$ is equivalent to an isobaric sum of unitary cuspidal automorphic representations; when $M/F$ is cyclic of prime degree $p$, for example, $\pi_M$ is either cuspidal or of the form $\boxplus_{j=0}^{p-1} (\eta\circ\tau^j)$, where $\tau$ is a generator of Gal$(M/F)$. Since base change is additive relative to isobaric sums, it follows that for any $\pi$ in $\mathcal A_u(n,F)$, $\pi_M$ lies in $\mathcal A_u(n,M)$.

\bigskip

\section{A Preliminary Step}

\medskip

\noindent{\bf Proposition 2.1} \, \it Let $F$ be a number field and $n\geq 1$ an integer. Suppose $\pi, \pi' \in \mathcal A_u(n,F)$ are such that for every positive integer $m\leq (n^2+1)/2$, and for all but a finite number of primes $v$ of $F$ of degree $m$, we have $\pi_v \simeq \pi'_v$. Then $\pi$ and $\pi'$ are isomorphic. \rm

\medskip

This is essentially an immediate consequence of the bound of Luo-Rudnick-Sarnak, as it has already been noted (and used) by Rajan for cuspidal representations in \cite{Raj2}. For completeness, we quickly go through the relevant points of \cite{Ra1} to make it evident that they carry over, modulo the basic results cited in section 1 and induction on the number of cuspidal isobaric summands, from ($n=2$; $\pi, \pi'$ cuspidal) to ($n$ arbitrary; $\pi, \pi'$ isobaric sums of unitary cuspidal representations).

\medskip

{\it Proof}. \, Denote by $X$ the complement in $\Sigma_F$ of the union of the archimedean places and the finite places where $\pi$ or $\pi'$ is ramified. Given any subset $Y$ of $X$ we set (as in \cite{Ra1}):
$$
Z_Y(s)=L_Y(\bar\pi\times\pi ,s)L_Y(\bar\pi'\times\pi',s)/L_Y(\bar\pi
\times\pi',s)L_Y(\bar\pi'\times\pi ,s).\leqno(2.1)
$$
Write
$$
\pi=\boxplus_{i=1}^\ell  m_i\pi_i, \, \, \, \pi'=\boxplus_{j=1}^r m_j' \pi_j',
$$
with $m_i, m_j' \in \N$, and $\pi_i$, $\pi'_j$ unitary cuspidal, with
$\pi_i\not\simeq \pi_a$ if $i\ne a$ and $\pi_j'\not\simeq \pi_b'$ if $j\ne b$.

\medskip

Suppose $\pi_i\not\simeq \pi'_j$ for all $i, j$.
Then, using Theorem 1.1,  we see that $Z_X(s)$ is holomorphic at every $s\neq 1$, with
$$
-{\rm ord}_{s=1}Z_X(s) \, = \, \mu+\mu',\leqno(2.2-a)
$$
where
$$
\mu=\sum_{i=1}^\ell m_i^2, \, \mu'=\sum_{j=1}^r {m'_j}^2.\leqno(2.2-b)
$$

We note that one knows (see \cite{HRa}) that $Z_Y(s)$ is of positive type, i.e., $\log Z_Y(s)$ is a Dirichlet series with non-negative coefficients.

\medskip

As the subproduct of an absolutely convergent Euler product is
absolutely convergent, we have the following

\medskip

\noindent{\bf Lemma 2.3} \, \it
Let $S$ denote the subset of $X$ consisting of finite places $v$ of degree $> \frac{n^2+1}{2}$.
Then the incomplete Euler products $L_S(\bar\pi\times\pi ,s)$ and $L_S
(\bar\pi\times\pi',s)L_s(\bar\pi'\times\pi ,s)$ converge absolutely
in $\{s\in\C \, | \, \Re(s)>1\}$.
\rm

\medskip

We may write
$$
\log (L_Y(\bar\pi\otimes\pi ,s))=\sum_{m\geqq 1}c_m(Y)m^{-s}\leqno(2.4)
$$
for all subsets $Y$ of $X$.  Then $c_m(Y)=0$ unless $m$ is of the form
$Nv^r$ for some $v\in Y$ and $r\in\N$, and when $m$ is of this
form,
$$
c_m(Y)=\sum_M \, \frac{1}{r}\sum_{1\leqq i,j\leqq 2}\overline{\alpha^r_{i,v}}
\alpha^r_{j,v}.
$$
where $M$ is the set of pairs $(v ,r)\in Y\times\N$ such that
$m=Nv^r$.

When $v \in S$, as $Nv > \frac{n^2+1}{2}$, the Luo-Rudnick-Sarnak bound (Theorem 1.2) implies that
$\sum_{m\geqq 1}c_m(S)m^{-s}$ converges in $\{\Re(s)\geq 1\}$.

One has a similar statement for $\log (L_S(\bar\pi'\otimes\pi ,s))$, $\log (L_S(\bar\pi'\otimes\pi ,s))$,
and $\log (L_S(\bar\pi'\otimes\pi' ,s))$. So we get the following

\medskip

\noindent{\bf Lemma 2.5} \it
Let $S$ be as in Lemma 2.3.  As $s$ goes to $1$ from the right on the real line, we have
$$
\log Z_S(s) \, = \, o\left(\log \frac{1}{s-1}\right).
$$
\rm

\medskip

Now, since $\pi_v \simeq \pi'_v$ for all but a finite number of places of $X$ outside $S$, we get, thanks to this Lemma, the following:
$$
\log Z_X(s) \, = \, 4\log L_X(\bar\pi\otimes\pi ,s) + o\left(\log \frac{1}{s-1}\right) \, = \,
4\log L_X(\bar\pi'\otimes\pi' ,s) + o\left(\log \frac{1}{s-1}\right).\leqno(2.6)
$$
Applying (2.2-b), we then get
$$
\mu \, = \, \mu',\leqno(2.7)
$$
and
$$
\log Z_X(s) \, = \, 4\mu\log \frac{1}{s-1} + o\left(\log \frac{1}{s-1}\right).\leqno(2.8)
$$
This contradicts (2.2-a) since $\mu=\mu' \geq 1$.

\medskip

Thus we must have
$\pi_i \simeq \pi'_j$ for {\it some} $(i,j)$. If $\pi$ or $\pi'$ is cuspidal, then both will need to be cuspidal with $\pi=\pi_i \simeq \pi'_j=\pi'$, and so we are done in this case. We may assume that $\pi, \pi'$ are non-cuspidal. Consider then the isobaric automorphic representations $\Pi$, $\Pi'$ such that
$$
\pi = \Pi\boxplus \pi_i, \, \pi' = \Pi'\boxplus \pi'_j.
$$
The $\Pi, \Pi'$ satisfy the hypotheses of Proposition 2.1, and we may find as before cuspidal isobaric summands $\pi_k$ of $\Pi$ and $\pi'_m$ of $\Pi'$ which are isomorphic. Continuing thus, by infinite decent, we arrive finally at the situation when one of the isobaric forms is cuspidal, which we have already taken care of. This proves Proposition 2.1.

\qed

\bigskip

\section{Central character and unitarity}

\medskip

Suppose $\pi$, $\pi'$ are cuspidal automorphic representations of GL$_n(\A_F)$ of respective central characters $\omega, \omega'$, such that $\pi_v \simeq \pi'_v$ for all but a finite number of primes $v$ of $F$ of degree $1$. Then $\omega$ and $\omega'$ agree at all (but a finite number of) the degree one places $v$, which forces the global identity
$$
\omega \, = \, \omega'.\leqno(3.1)
$$
In fact, by Hecke, this conclusion will result as soon as $\omega$ and $\omega'$ agree at a set of primes of density $> 1/2$.

\medskip

It is a standard fact that, given a cuspidal $\pi$, there is a unique real number $t(\pi)$ such that $\pi\otimes\vert\cdot\vert^{-t(\pi)}$ is unitary; here $\vert\cdot\vert$ denotes the $1$-dimensional representation $g\mapsto \vert{\rm det}(g)\vert$. Taking central characters, we see then that $\omega\vert\cdot\vert^{-nt(\pi)}$ is a unitary character. Thanks to (3.1), we will then get
$$
t(\pi) \, = \, t(\pi').\leqno(3.2)
$$

This allows us, in the proof of Theorem A, to assume that $\pi, \pi'$ are unitary cuspidal automorphic representations.

\bigskip

\section{Nested chains of cyclic $p^2$-extensions}

\medskip

Let $p$ be a prime. We will call an extension $L/F$ of number fields of degree $p^r$, for some $r \geq 2$, a {\it nested chain of cyclic $p^2$-extensions} if there is an increasing filtration of fields
$$
F=L_0 \subset L_1 \subset L_2 \subset \dots \subset L_{r-2} \subset L_{r-1} \subset L_r=L,\leqno(4.1)
$$
with
$$
[L_j : L_{j-1}] = p, \, \, \forall \, j \in\{1, 2, \dots, r\},\leqno(4.2)
$$
and
$$
L_j/L_{j-2}: \, \, {\rm cyclic}, \, \, \forall \, j \in\{2, \dots, r\}.\leqno(4.3)
$$

\medskip

An easy example is given by a cyclic $p^r$ extension, while a better example is the following. Let $F$ contain $\mu_{p^2}$. (As usual, we write $\mu_n$ for the group of $n$-th roots of unity in the algebraic closure of $F$.) Let $\alpha$ be an element of $F$ which is not a $p$-th power. Put $\alpha_0=\alpha$ and define $\alpha_j$, for $j=1, \dots, r$, recursively by taking it to be a $p$-th root of $\alpha_{j-1}$, and set $L_j=L_{j-1}(\alpha_j)$ and $L_0=F$. Note that for $j \geq 2$, $L_j/L_{j-2}$ is cyclic of order $p^2$ by Kummer theory, because $\alpha_j^{p^2}=\alpha_{j-2}$, and $\mu_{p^2}\subset L_{j-2}$, making all the conjugates of $\alpha_j$ over $L_{j-2}$ lie in $L_j$. (For this example, it is in fact sufficient to have $\mu_p \subset F$ and $\mu_{p^2}\subset L_1$, as seen by the case $L_1=F(\mu_{p^2})$.)

\medskip

\noindent{\bf Lemma 4.4} \, \it Let $L/F$ be a nested chain of cyclic $p^2$-extensions (of number fields), with $[L:F]=p^r$ and filtration $\{L_j\}$ as above. Suppose $v_0$ is a finite place of $F$, unramified in $L$, which is inert in $L_1$. Then there exists, for each $j\geq 1$, a unique place $v_j$ of $L_j$ lying over $v_{j-1}$, so that $Nv_j=(Nv_{j-1})^p$. In particular, $Nv_r=(Nv_0)^{p^r}$.\rm

\medskip

{\it Proof}. \, Let us first treat the case when $r=2$, i.e., when $L/F$ is cyclic of degree $p^2$. Since $v_0$ is inert in the intermediate field $L_1$, we need to check that $v_0$ does not split into $p$ places in $L$. Suppose, to the contrary, that it does split that way. Let $u$ be one of the $p$ places of $L$ above $v_0$. It must then be fixed by a subgroup $H$ of Gal$(L/F)$ of order $p$, with $H$ giving the local Galois group ${\rm Gal}(L_u/F_{v_0})$. Since $v_0$ is inert in $L_1$ with divisor $v_1$, $u$ necessarily has degree $1$ over $v_1$, and so $H = {\rm Gal}(L_{1, v_1}/F_{v_0})$. If $\sigma$ is a non-trivial element of $H$, then it acts non-trivially on $L_{1,v_1}$, and hence on $L_1$. On the other hand, since $L/F$ is cyclic, it has a unique subgroup of order $p$, which forces $H$ to be Gal$(L/L_1)$, implying that $\sigma$ acts trivially on $L_1$, yielding a contradiction. Put another way, if $v_0$ has degree $p$ in $L$, then the corresponding Frobenius class $Fr_{v_0}$ is given by an element $\sigma$ of Gal$(L/F)$ of order $p$, which has trivial image in the quotient by $H=\langle \sigma\rangle$, making $v_0$ split in the fixed field $L^H$ of $H$. Clearly, $L^H$ must be $L_1$ by the cyclicity of $L/F$. Either way, the case $r=2$ is now settled.

Now let $r > 2$, and assume by induction that the Lemma holds for $r-1$. So for every $j \leq r-1$, there is a unique place $v_j$ of $L_j$ above $v_{j-1}$ (of $L_{j-1}$). Now all we have to show is that $v_{r-1}$ is inert in $L=L_r$. Since $L_r/L_{r-2}$ is cyclic of order $p^2$, and since (by induction) the place $v_{r-2}$ of $L_{r-2}$ is inert in $L_{r-1}$, we conclude what we want by appealing again to the $r=2$ scenario.

The assertion about the norm of $v_r$ follows.

\qed

\medskip

\noindent{\bf Lemma 4.5} \, \it Let $L^{(i)}/F$, $1\leq i \leq k$ be disjoint $p^r$-extensions. Suppose moreover that every $L^{(i)}$ is a nested chain of cyclic $p^2$-extensions with respective filtrations
$$
F=L_0^{(i)}\subset L_1^{(i)} \subset \dots \subset L_r^{(i)}=L^{(i)}.
$$
Let $v_0^{(i)}$, $1\leq i \leq k$, be distinct primes of $F$, unramified in the compositum $M:=L^{(1)}L^{(2)}\dots L^{(k)}$, such that each $v_0^{(i)}$ is inert in $L_1^{(i)}$. Then, if ${\tilde v}^{(i)}$ is a prime of $M$ lying above $v_0^{(i)}$, we have
$$
N{\tilde v}^{(i)} \, \geq \, (Nv_0^{(i)})^{p^r}, \, \, \, \forall \, i\leq k.
$$
\rm

\medskip

{\it Proof}. \, Fix any $i \leq k$. By Lemma 4.4, for each $j\geq 2$, there is a unique prime $v_j^{(i)}$, of $L_j^{(i)}$ lying above $v_{j-1}^{(i)}$. Then $\tilde v^{(i)}$ must lie above $v_r^{(i)}$ in the extension $M/L^{(i)}$. So
$$
N{\tilde v^{(i)}} \, \geq \, Nv_r^{(i)}.\leqno(4.6)
$$
On the other hand, by Lemma 4.4, we have
$$
Nv_r^{(i)} \, = \, (Nv_0^{(i)})^{p^r}.\leqno(4.7)
$$
The assertion of Lemma 4.5 now follows by combining (4.6) and (4.7).

\qed

\bigskip

\section{Isomorphism over suitable solvable extensions $L/K$, $L\supset E$}

\medskip

Let $K/k$ be a cyclic $p$-extension. For $j\geq 1$, denote by $\Sigma^j_{K/k}$ the set of finite places $v$ of $K$ which are unramified over $k$ and of degree $j$ over $k$; of course this set is non-empty only for $j\in \{1,p\}$. Let $\pi, \pi'$ be cuspidal automorphic representations of GL$_n(\A_K)$ such that, as in the setup of Theorem A,
$$
\pi_v \, \simeq \, \pi'_v, \, \, \foralmostall \, v \in \Sigma^1_{K/k}.\leqno(5.1)
$$

As noted in section 3, the central characters of $\pi$ and $\pi'$ must be the same, and moreover, we may assume that $\pi,  \pi'$ are unitary.

\medskip

If $p > (n^2+1)/2$, then Theorem A follows immediately from Proposition 2.1. In general, fix a positive integer $r$ such that
$$
p^r \,  > \, (n^2+1)/2.\leqno(5.2)
$$

The object of this section is to prove the following:

\medskip

\noindent{\bf Proposition 5.3} \, \it Let $K/k$, $\pi, \pi'$ be as in Theorem A. Then there is a finite solvable extension $L/K$ containing $E:=K(\mu_{p^2})$ such that the base changes $\pi_L$, $\pi'_L$, satisfy
$$
\pi_L \, \simeq \, \pi'_L.
$$
\rm

\medskip

In fact the number field $L$ we construct below will be much nicer than just being solvable over $K$. The extension $L/E$ will turn out to be the compositum of a finite number of of $p^r$-extensions $L^{(i)}$ , with each of them a nested chain of cyclic $p^2$-extensions. The Galois closure of $L$ over $K(\mu_p)$ will again be a $p$-power extension, hence nilpotent. We will also have some freedom in the choice of the $L^{(i)}$, and their filtrations, which will become relevant in the next section when we descend to $E$.

\medskip

Put $K'=K(\mu_p)$ and $k'=k(\mu_p)$. Then $K'/k'$ is still a cyclic $p$-extension. The following Lemma is clear since $K'/K$ and $k'/k$ are of degree dividing $p-1$.

\medskip

\noindent{\bf Lemma 5.4} \, \it Let $v \in \Sigma^j_{K/k}$, for $j\in\{1,p\}$. Then, for every prime $v'$ of $K'$ above $v$, we have $v' \in \Sigma^j_{K'/k'}$.\rm

\medskip

Consequently, the hypotheses of Theorem A are preserved for $K'/k'$, and we may assume from here on, after replacing $k$ (resp. $K$) by $k'$ (resp. $K'$), that
$$
\mu_p \, \subset \, k.\leqno(5.5)
$$

\medskip

{\bf Proof of Proposition 5.3 when $K=E$}

\medskip

Since $\mu_p \subset k$, we may realize the cyclic $p$-extension $K$ as $k(\alpha^{1/p})$, for an element $\alpha$ in $k$ which is not a $p$-th power (in $k$). Now fix a positive integer $r$ for which (5.2) holds. Choose a sequence of elements $\alpha_{-1}=\alpha$, $\alpha_0, \dots, \alpha_r$ in the algebraic closure of $K$, and the corresponding chain of fields $k=L_{-1}, K=L_0, \dots, L_r$ such that for each $j\geq 0$,
$$
L_j=L_{j-1}(\alpha_j), \, \, {\rm with} \, \, \alpha_j^p = \alpha_{j-1}.\leqno(5.6)
$$
Clearly, every $L_j/L_{j-1}$ is cyclic of order $p$, and so $[L_r:K]=p^r$. Moreover, since $\mu_{p^2}\subset E=K$, each $L_j/L_{j-2}$ is also cyclic by Kummer theory. In other words, $L_r/K$ is a nested chain of cyclic $p^2$-extensions. In fact, $L_r/k$ is also such a nested chain, but of degree $p^{r+1}$.

Now put $L=L_r$. Applying Lemma 4.4, we then see that for every prime $\tilde v$ in $L$ lying over some $v$ in $\Sigma^p_{K/k}$, the degree of $\tilde v$ is $p^r$ over $k$, hence has degree at least $p^r$ over $\Q$. On the other hand, every other prime $\tilde u$ of $L$ unramified over $k$ lies above some $u$ in $\Sigma^1_{K/k}$. So the hypotheses of Theorem A imply (by base change \cite{AC}) that $\pi_{L,\tilde u}\simeq \pi'_{L,\tilde u}$. (Such a $\tilde u$ could have small degree, like $p$, over $K$, but nevertheless it must lie over a prime $u$ of degree $1$ over $k$, which is all that matters to us.) Putting these together, and applying Proposition 2.1 over $L$, we get Proposition 5.3 when $K=E$.
\qed

\medskip

{\bf Proof of Proposition 5.3 when $K\ne E$}

\medskip

Here we want to base change and consider the cyclic $p$-extension
$$
E/F, \, \, {\rm with} \, \, F = k(\mu_{p^2}), \, E=KF.\leqno(5.7)
$$
Clearly, the $(p,p)$-extension $E/k$ contains $p+1$ subfields $F^{(i)}$, $0\leq i \leq p$, of degree $p$ over $k$, with one of them being $K$; say $K=F^{(0)}$. We need the following

\medskip

\noindent{\bf Lemma 5.8} \, \it  Let $v\in \Sigma^p_{K/k}$ be unramified in $E$. Then $v$ splits into $p$ places $v_1, \dots, v_p$ in $E$, and there is a (unique) cyclic $p$-extension $F^{(i)}$ of $k$ (depending on $v$), $1\leq i \leq p$, such that each $v_j$ lies in $\Sigma^p_{E/F^{(i)}}$. In other words, if $z$ is the unique place of $k$ below $v$, then $z$ splits into $p$ places in $F^{(i)}$, each of which is inert in $E$.
\rm

\medskip

{\it Proof of Lemma 5.8}. \, Since $G:=$Gal$(E/k)$ is $\Z/p \times \Z/p$, the decomposition groups are either trivial or of order $p$. So, if $z$ is the place of $k$ lying below $v$, its Frobenius class $Fr_z$ in $G$ is given by an element $\sigma$ of order $p$ (since $z$ is inert in $K$). So $v$ must split in $K$. If we put $H=\langle \sigma\rangle$, then $K^H$ is $F^{(i)}$ for a unique $i\in \{1, \dots, p\}$. Then $z$ splits in $F^{(i)}$ and then becomes inert in $E$, as claimed.
\qed

\medskip

Fix an index $i\in \{1, \dots, p\}$. As $\mu_p\subset k\subset F^{(i)}$, we may find an element $\alpha^{(i)}$ in $F^{(i)}$ which is not a $p$-th power such that
$$
E \, = \, F^{(i)}((\alpha^{(i)})^{1/p}).\leqno(5.9)
$$
Choose a sequence of elements $\alpha_{-1}^{(i)}=\alpha^{(i)}$, $\alpha_0^{(i)}, \dots, \alpha_r^{(i)}$ in the algebraic closure of $E$, and the corresponding chain of fields $F^{(i)}=L_{-1}^{(i)}, E=L_0^{(i)}, \dots, L_r^{(i)}$ such that for each $j\geq 0$,
$$
L_j^{(i)}=L_{j-1}^{(i)}(\alpha_j^{(i)}), \, \, {\rm with} \, \, (\alpha_j^{(i)})^p = \alpha_{j-1}^{(i)}.\leqno(5.10)
$$
By construction, every $L_j^{(i)}/L_{j-1}^{(i)}$ is cyclic of order $p$, and so $[L_r^{(i)}:E]=p^r$. Moreover, since $\mu_{p^2}\subset E$, each $L_j^{(i)}/L_{j-2}^{(i)}$ is also cyclic by Kummer theory. In other words, $L_r^{(i)}/E$ is a nested chain of cyclic $p^2$-extensions. In fact, $L_r^{(i)}/F^{(i)}$ is also such a nested chain (of degree $p^{r+1}$).

This way we get $p$ nested chains $L^{(i)}/E$, disjoint over $K$ from each other. Let $L$ be the compositum of the $L^{(i)}$, as $i$ runs over $\{1, \dots, p\}$. Pick any place $v$ in $\Sigma^p_{K/k}$. Then we know (by Lemma 5.8) that there is a unique $i \leq p$ such that each of the divisors $v_k$ of $v$ in $E$, $1\leq k \leq p$, lies in $\Sigma^p_{E/L^{(i)}}$. Then by the $r=2$ case of Lemma 4.4, $v_k$ is inert in $L^{(1)}$. Applying Lemma 4.5, we then see that every prime $\tilde v$ of $L$ lying over some $v_k$ (and hence over $v$) is of degree $\geq p^r > (n^2+1)/2$.  So one may apply Proposition 2.1 and conclude that $\pi_L$ and $\pi'_L$ are isomorphic.

\qed

\bigskip

\section{Descent to $E=K(\mu_{p^2})$}

\medskip

Let us preserve the notations of the previous section. Thanks to Proposition 5.3, we know that for the $p$-power extension $L/E$ we constructed there, one has
$$
\pi_L \, \simeq \pi'_L.\leqno(6.1)
$$
In order to prove Theorem A, we need to descend this isomorphism down to $E$. For this we will make use of the fact that there is quite a bit of freedom in choosing $L$.

\medskip

{\bf Proof of descent when $K=E$}

\medskip

After realizing $E$ as $k(\alpha^{1/p})$ for some $\alpha$ ($=\alpha_{-1}$) in $k$ which is not a $p$-th power, we chose a sequence of elements $\alpha_j, 0\leq j \leq r$, with $\alpha_j=\alpha_{j-1}^{1/p}$, and set $L_j=L_{j-1}(\alpha_j)$. We may replace $\alpha$ by $\alpha\beta^p$ for any $\beta$ in $k-k^{p}$, which will have the effect of leaving $E=L_0$ intact, but changing $L_1$ from $E(\alpha_1)$ to $E(\alpha_1\beta_1)$ for a $p$-th root $\beta_1$ of $\beta$. Using this we can ensure, for a suitable choice of $\beta$, that the discriminant of $L_1/E$ is divisible by a prime $P_1$ not dividing the conductor of either $\pi_E$ or $\pi'_E$. Next we may choose a $\gamma \in k-k^p$ and put $\alpha_0=\alpha_0\beta^p\gamma^{p^2}$, which will not change $L_0$ and $L_1$, but will change $L_2$, and we may arrange for the discriminant of the new $L_2/L_1$ to be divisible by a prime $P_2$ of $L_1$ whose norm down to $E$ is relatively prime to ${\mathfrak c}(\pi_E){\mathfrak c}(\pi'_
 E)P_1$. This way we may continue and modify all the $L_j$ so that at each stage $L_j/L_{j-1}$, the relative discriminant is divisible by a new prime $P_j$ of $L_{j-1}$ whose norm down to $E$ is relatively prime to ${\mathfrak c}(\pi_E){\mathfrak c}(\pi'_E)P_1N_{L_1/E}(P_2)\dots N_{L_{j-2}/E}(P_{j-1})$.

Now look at the top stage $L_r/L_{r-1}$. Thanks to (6.1), we know by the properties of base change (\cite{AC}) that every cuspidal isobaric component $\eta$, say, of $\pi_{L_{r-1}}$ will be twist equivalent to a cuspidal isobaric component $\eta'$ of $\pi'_{L_{r-1}}$. More precisely, we will need to have, for some integer $j$ mod $p$,
$$
\eta' \, \simeq \, \eta\otimes\delta_r^j,\leqno(6.2)
$$
where $\delta_r$ is the character of order $p$ of (the idele classes of) $L_{r-1}$ attached to $L_r$. But the conductor of $\delta_r$ is divisible by the prime $P_r$, whose norm down to $E$ is, by construction, relatively prime to the conductors of $\pi_E$ and $\pi'_E$ and to the discriminant of $L_{r-1}/E$. This forces $j=0$, i.e., $\eta \simeq \eta'$. Peeling off this way isomorphic cuspidal components of $\pi_{L_{r-1}}$ and $\pi'_{L_{r-1}}$ one at a time, we conclude that $\pi_{L_{r-1}}$ is isomorphic to $\pi'_{L_{r-1}}$. Next, by induction on $r-j$, we deduce similarly that, for every $j\in \{0, \dots, r-1\}$,
$$
\pi_{L_j} \, \simeq \, \pi'_{L_j},\leqno(6.3)
$$
which proves the assertion of Theorem A.
\qed

\medskip

{\bf Proof of descent when $K\ne E$}

\medskip

For each $i=\{1, \dots, p\}$, we may modify the elements $\alpha_j^{(i)}$ and thus the fields $L_j^{(i)}$ as above, with a new prime divisor $P_j^{(i)}$ of the discriminant of $L_j/L_{j-1}$ popping up at stage $j$, which is prime to the conductors of $\pi_E$, $\pi'_E$, and the discriminant of $L_{j-1}/E$. Now we may, and we will, also choose these primes in such a way that the sets $\{P_1^{(i)}, \dots, P_r^{(i)}\}$ and $\{P_1^{(k)}, \dots, P_r^{(k)}\}$ are disjoint whenever $i \ne k$. Now we may realize $L$ as a sequence of cyclic $p$-extensions, such that at each stage there is a new prime divisor of the relative discriminant. We may then descend each step as above and finally conclude that
$$
\pi_E \, \simeq \, \pi'_E,\leqno(6.4)
$$
as asserted.
\qed

\bigskip

\section{Descent to $K(\mu_p)$}

\medskip

As before, we may assume that $\mu_p \subset k\subset K$. If $\mu_{p^2}\subset K$, i.e., if $E=K$, then we have already seen above that we have an isomorphism $\pi\simeq \pi'$ over $K$.

So we may, and we will, assume below that $K\ne E$. Then
$$
E=KF, \, k=K\cap F, \, \, \, {\rm where} \, \, F=k(\mu_{p^2}),\leqno(7.1)
$$
with
$$
[E:F]=[K:k]=[E:K]=[F:k]=p,
$$
and by section 6,
$$
\pi_E \, \simeq \, \pi'_E.\leqno(7.2)
$$
This implies that if $v$ is any prime of $K$ which splits into $p$ primes $w_1, \dots, w_p$ in $E$, then by \cite{AC}, we have ($\forall j\leq p$)
$$
\pi_v \, \simeq \, \pi_{w_j} \, \simeq \, \pi'_{w_j} \, \simeq \, \pi'_{v}.\leqno(7.3)
$$

On the other hand, since $E/k$ is a $(p,p)$-extension, in particular not cyclic of order $p^2$, any prime $u$ of $k$ which is inert in $K$ must split in $E$ (assuming $u$ is unramified in $E$). This implies, thanks to (7.3), the following:
$$
\pi_v \, \simeq \, \pi'_v, \, \, \, \forall \, v\in \Sigma^p_{K/k}-{{\rm finite} \, \, \, {\rm set}}.\leqno(7.4)
$$

When we combine (7.4) with the hypothesis of Theorem A that
$$
\pi_v \, \simeq \, \pi'_v, \, \, \forall \, v \in \Sigma^1_{K/k},\leqno(7.5)
$$
we immediately get the desired isomorphism
$$
\pi \, \simeq \, \pi' \, \, ({\rm over} \, \, K).
$$

We are now done with the proof of Theorem A. The assertion of Corollary B is obvious given Theorem A (since $\mu_2\subset \Q\subset K$).

\qed

\vskip 0.2in

Dinakar Ramakrishnan

253-37 Caltech

Pasadena, CA 91125, USA.

dinakar@caltech.edu

\bigskip

\end{document}